\documentclass{amsart}
\usepackage{amssymb}
\usepackage{amsfonts}
\usepackage{latexsym}

\newtheorem{theorem}{Theorem}[section]
\newtheorem{lemma}[theorem]{Lemma}
\newtheorem{proposition}[theorem]{Proposition}
\newtheorem{corollary}[theorem]{Corollary} 
\theoremstyle{definition}  
\newtheorem{definition}[theorem]{Definition}
\newtheorem{example}[theorem]{Example}
\newtheorem{conjecture}[theorem]{Conjecture}  
\newtheorem{remark}[theorem]{Remark}
\newtheorem{note}[theorem]{Note}

\newcommand{\Tr}{\text{Tr}}
\newcommand{\id}{\text{id}}

\newcommand{\End}{\text{End}} 

\newcommand{\Hom}{\text{Hom}} 
\newcommand{\Ad}{\text{Ad}}

\newcommand{\Rep}{\text{Rep}}

\newcommand{\op}{\text{op}}

\newcommand{\eps}{\varepsilon}

\newcommand{\I}{\mathcal{I}}

\newcommand{\G}{\widetilde{G}}

\newcommand{\ov}{\overline}

\newcommand{\actl}{\rightharpoonup}
\newcommand{\actr}{\leftharpoonup}

\newcommand{\la}{\langle\,} 
\newcommand{\ra}{\,\rangle}

\newcommand{\1}{_{(1)}} 
\newcommand{\2}{_{(2)}} 
\newcommand{\3}{_{(3)}} 
 
\renewcommand{\I}{^{(1)}} 
\newcommand{\II}{^{(2)}}

\newcommand{\bTheta}{{\bar \Theta}}

\begin{document}

\title{On the structure of weak Hopf algebras}

\author{Dmitri Nikshych} 
\address{MIT, Department of Mathematics,
77 Massachusetts Avenue, Room 2-130, Cambridge, MA 02139-4307 USA}
\email{nikshych@math.mit.edu}
\thanks{The author is grateful to Pavel Etingof for 
numerous stimulating discussions, advice,  and support. 
This research was partially conducted by the author for the
Clay Mathematics Institute.
The author was  partially supported by the NSF grant  DMS-9988796 and
thanks MIT for the hospitality.}
  
\date{June 1, 2001}
\begin{abstract}
We study the group of group-like elements of a weak Hopf algebra
and derive an analogue of Radford's formula for the fourth power
of the antipode $S$, which implies that the antipode has a finite
order modulo a trivial automorphism.  We find a sufficient condition
in terms of $\Tr(S^2)$ for a weak Hopf algebra to be semisimple, 
discuss relation between semisimplicity and cosemisimplicity, and
apply our results to show that a dynamical twisting deformation
of a semisimple Hopf algebra is cosemisimple.
\end{abstract} 
\maketitle  


\begin{section}
{Introduction}

Weak Hopf algebras were introduced by G.~B\"ohm and K.~Szlach\'anyi in
\cite{BSz1} and \cite{Sz} (see also their joint work \cite{BNSz}
with F.~Nill) as a generalization of ordinary Hopf algebras and groupoid 
algebras. A weak Hopf algebra is a vector space that has both algebra 
and coalgebra structures related to each other in a certain self-dual way
and that possesses an antipode.  The main difference
between ordinary and weak Hopf algebras comes from the fact that
the comultiplication of the latter is no longer required to preserve
the unit (equivalently, the counit is not required to
be an algebra  homomorphism) and results in the existence of two canonical
subalgebras playing the role of ``non-commutative bases'' in
a ``quantum groupoid''. The axioms of a weak Hopf 
algebra are self-dual, which ensures that 
when $H$ is finite dimensional, the dual vector space $H^*$
has a natural structure of a weak Hopf algebra.

For the foundations of the weak Hopf algebra theory we refer the
reader to \cite{BNSz} and to the recent survey \cite{NV4}.

The initial motivation to study weak Hopf algebras was their connection
with the theory of algebra extensions. It was explained in
\cite{NV2,NV3} that  weak Hopf $C^*$-algebras naturally arise as symmetries
of finite depth von Neumann subfactors. A purely algebraic analogue
of this result was proved in \cite{KN}, where it
was shown that a depth two Frobenius extension  of algebras $A\subset B$
with a separable centralizer $C_B(A)$ comes from a smash product with
a semisimple and cosemisimple weak Hopf algebra $H$, i.e., 
$(A\subset B) \cong (A \subset A\# H)$. 

Another important application of weak Hopf algebras is that they
provide a natural framework for the study of dynamical twists \cite{ES} in
Hopf algebras. It was proved in  \cite{EN1} that every dynamical twist in
a Hopf algebra gives rise to a weak Hopf algebra. Also in  \cite{EN1}
a family of dynamical quantum groups (weak Hopf algebras
corresponding to dynamical twists in quantum groups at roots of unity)
was constructed. These weak Hopf algebras were shown to be quasitriangular
with non-degenerate $R$-matrices and, therefore, self-dual. 

It turns out that many important properties of ordinary Hopf algebras
have ``weak'' analogues. For example, the category $\Rep(H)$ of finite rank
left modules over a weak Hopf algebra $H$ is a rigid monoidal category
(the category $\Rep(H)$ was defined and studied in \cite{BSz2} for
weak Hopf $C^*$-algebras and in \cite{NTV} for general weak Hopf algebras). 
The importance of weak Hopf algebra representation categories
can be seen from the result of V.~Ostrik \cite{O},
who proved that {\em every} semisimple rigid monoidal category with finitely
many (classes of) simple objects is equivalent to $\Rep(H)$ for some 
semisimple weak Hopf algebra $H$.

Also, the theory of integrals for weak Hopf algebras developed in \cite{BNSz}
is essentially parallel to that of ordinary Hopf algebras. Using it,
one can prove an analogue of Maschke's theorem \cite[Theorem 3.13]{BNSz}
for weak Hopf algebras and show that semisimple weak Hopf algebras are
finite dimensional.

Despite these similarities, the structure of weak Hopf algebras is
much more complicated than that of usual Hopf algebras, 
even in the semisimple case. For instance, the antipode of 
a semisimple weak Hopf algebra $H$ over $\mathbb{C}$,
the field of complex numbers, may have an infinite order and 
quantum dimensions of  irreducible $H$-modules can be non-integer 
(this is clear from the relation between  weak Hopf algebras and semisimple 
monoidal categories mentioned above). 
Also, weak Hopf algebras of prime dimension can be non-commutative 
and non-cocommutative, quite in contrast with the usual Hopf algebra theory, 
cf.\ \cite{Z}. An example of an indecomposable
semisimple weak Hopf algebra of dimension 
$13$ having such properties was given in \cite{BSz1}
(see also \cite[Appendix]{NV3}).

In this paper we start an investigation of the structure of finite dimensional
weak Hopf algebras and prove weak Hopf algebra analogues of several classical
Hopf algebra results obtained in \cite{LR2} and \cite{R1}.

After giving necessary definitions and discussing basic properties
of weak Hopf algebras in Section $2$, we classify all {\em  minimal} 
weak Hopf algebras (i.e., those generated by the identity element)
in Section $3$. We show that every weak Hopf algebra can be obtained as
a deformation of a weak Hopf algebra with the property that the square
of the antipode acts as the identity on the minimal weak Hopf subalgebra.

Next, in Section $4$ we define the group $G(H)$ of group-like elements 
of a weak Hopf algebra $H$ and show that it contains a normal subgroup 
$G_0(H)$ of {\em trivial} group-like elements that belong to the minimal
weak Hopf subalgebra of $H$.  An adjoint action
of any group-like element defines a weak Hopf algebra  
automorphism of $H$. Note that even in finite dimensional case both 
$G(H)$ and $G_0(H)$ are usually infinite.
However, when $H$ is finite dimensional, we show that the quotient group 
$\G(H) = G(H)/ G_0(H)$ is finite. This $\G(H)$ turns out to be the 
correct weak Hopf algebra
analogue of the group of group-like elements in an ordinary Hopf algebra.
We introduce the notion of a  distinguished coset of group-like elements of 
$H$, that ``measures the difference'' between left and right 
integrals in $H^*$. 

In Section $5$ we extend to weak Hopf algebras the result of 
D.~Radford \cite{R1} stating that the antipode of a finite dimensional 
Hopf algebra has a finite order.
Namely, we establish a formula for the fourth power of the antipode
analogous to \cite[Proposition 6]{R1} and 
show that the order of the antipode of a   
finite dimensional Frobenius weak Hopf algebra is finite modulo the group of 
trivial automorphisms of $H$. As in \cite{R1}, the proof uses 
non-degenerate integrals and distinguished  group-like elements 
of $H$ and $H^*$. 

Finally, we extend the Larson-Radford formula for $\Tr(S^2)$ 
\cite{LR2} and use it to establish a sufficient condition for a weak Hopf 
algebra being both semisimple and cosemisimple. We also 
prove that semisimplicity of a weak Hopf algebra with coinciding bases
is equivalent to its cosemisimplicity. As an application, we show that 
a dynamical twisting deformation of a semisimple Hopf algebra \cite{EN1}
is a semisimple and cosemisimple weak Hopf algebra (dynamical twists
are closely related to the dynamical quantum Yang-Baxter equation
of Felder \cite{F}). 

\end{section}


\begin{section}
{Preliminaries}
\label{Prelim}


Throughout this paper $k$ denotes a field. We use Sweedler's notation
for a comultiplication:  $\Delta(c) = c\1 \otimes c\2$. For an algebra
$A$ we denote by $Z(A)$ its center.

\subsection*{Definition of a weak Hopf algebra}
Below we collect the definition and basic properties of
weak Hopf algebras.

\begin{definition}[\cite{BNSz}, \cite{BSz1}] 
\label{finite weak Hopf algebra}
A {\em weak Hopf algebra} is a vector space $H$
with the structures of an associative algebra $(H,\,m,\,1)$ 
with a multiplication $m:H\otimes_k H\to H$ and unit $1\in H$ and a 
coassociative coalgebra $(H,\,\Delta,\,\epsilon)$ with a comultiplication
$\Delta:H\to H\otimes_k H$ and counit $\epsilon:H\to k$ such that:
\begin{enumerate}
\item[(i)] The comultiplication $\Delta$ 
is a (not necessarily unit-preserving) homomorphism of algebras:
\begin{equation}
\label{Delta m}
\Delta(hg) = \Delta(h) \Delta(g), \qquad h,g\in H,
\end{equation}
\item[(ii)] The  unit and counit satisfy the following identities: 
\begin{eqnarray}
\label{Delta 1}
(\Delta \otimes \id) \Delta(1) 
& =& (\Delta(1)\otimes 1)(1\otimes \Delta(1)) 
= (1\otimes \Delta(1))(\Delta(1)\otimes 1), \\
\label{eps m}
\epsilon(fgh) &=& \epsilon(fg\1)\, \epsilon(g\2h) = \epsilon(fg\2)\, 
\epsilon(g\1h),
\end{eqnarray}
for all $f,g,h\in H$. 
\item[(iii)]
There is a linear map $S: H \to H$, called an {\em antipode}, such that
\begin{eqnarray}
m(\id \otimes S)\Delta(h) &=&(\epsilon\otimes\id)(\Delta(1)(h\otimes 1)),
\label{S epst} \\
m(S\otimes \id)\Delta(h) &=& (\id \otimes \epsilon)((1\otimes h)\Delta(1)),
\label{S epss} \\
S(h) &=& S(h\1)h\2 S(h\3),
\label{S id S}
\end{eqnarray}
for all $h\in H$.
\end{enumerate} 
\end{definition}

Axioms (\ref{Delta 1})  and (\ref{eps m}) above 
are  analogous to the usual bialgebra axioms of 
$\Delta$ being a unit preserving map and $\epsilon$ being an
algebra homomorphism. Axioms (\ref{S epst}) and (\ref{S epss}) 
generalize the properties
of the antipode with respect to the counit. Also, it is possible to show
that given (\ref{Delta m}) - (\ref{S epss}),
axiom (\ref{S id S}) is equivalent to $S$ being both
anti-algebra and anti-coalgebra map.

\begin{remark}
A weak Hopf algebra is a Hopf algebra if and only if the comultiplication 
is unit-preserving and if and only if $\epsilon$ is a homomorphism of algebras.
\end{remark}

A {\em morphism} between weak Hopf algebras $H_1$ and $H_2$
is a map $\alpha : H_1 \to H_2$ which is both algebra and coalgebra
homomorphism preserving $1$ and $\epsilon$
and which intertwines the antipodes of $H_1$ and $H_2$,
i.e., $ \alpha\circ S_1 = S_2\circ \alpha$. The image of a morphism is
clearly a weak Hopf algebra. 

When $\dim_k H < \infty$, there is a natural weak Hopf algebra
structure on the dual vector space $H^*=\Hom_k(H,k)$ given by
\begin{eqnarray}
& & \la \phi\psi,\,h \ra = \la \phi\otimes\psi,\, \Delta(h)\ra, \\
& & \la \Delta(\phi),\,h\otimes g \ra =  \la \phi,\,hg \ra, \\
& & \la S(\phi),\,h \ra = \la \phi,\,S(h)\ra,
\end{eqnarray}
for all $\phi,\psi \in H^*,\, h,g\in H$. The unit
of $H^*$ is $\epsilon$  and the counit is $\phi \mapsto \la\phi,\, 1\ra$.

In what follows we use the Sweedler arrows for the dual actions,
writing 
\begin{equation}
\label{Sweedler arrows}
h\actl\phi = \phi\1 \la \phi\2,\, h\ra,
\qquad
\phi\actr h =\la \phi\1,\,h \ra \phi\2.
\end{equation}   
for all $h\in H,\phi\in H^*$.

\begin{example}
\label{groupoid examples}
Let $G$ be a {\em groupoid} (a small category with inverses)
with finitely many objects, then the groupoid
algebra $kG$ (generated by morphisms $g\in G$ with the product of 
two morphisms being equal to  their composition
if the latter is defined and $0$ otherwise) 
is a weak Hopf algebra via :
\begin{equation}
\Delta(g) = g\otimes g,\quad \epsilon(g) =1,\quad S(g)=g^{-1},
\quad \mbox{ for all } g\in G.
\label{groupoid algebra}
\end{equation}
In fact, when $G$ is finite and $k$ is algebraically closed and has
characteristic $0$, this is the most general example of
a cocommutative weak Hopf algebra \cite{thesis}.

If $G$ is finite, then  the dual weak Hopf algebra $(kG)^*$ is
isomorphic to the algebra of functions on $G$, i.e.,
it is generated by idempotents  $p_g,\, g\in G$ such that
$p_g p_h= \delta_{g,h}p_g$, with the following structure operations
\begin{equation}
\Delta(p_g) =\sum_{uv=g}\,p_u\otimes p_v,\quad \epsilon(p_g)
            = \delta_{g,gg^{-1}}, \quad S(p_g) =p_{g^{-1}}.
\label{dual groupoid algebra}
\end{equation}
This is the most general example of a finite dimensional commutative
weak Hopf algebra over an algebraically closed field of
characteristic $0$ \cite{thesis}.

More interesting examples of weak Hopf algebras arise as dynamical twisting
deformations of Hopf algebras \cite{EN1} and as symmetries of 
finite depth subfactors \cite{NV2}, \cite{NV3}.
\end{example}

We refer the reader to the original article \cite{BNSz} and recent 
survey \cite{NV4} for a detailed introduction to the theory of
weak Hopf algebras and its applications.

\subsection*{Counital maps and bases}

The linear maps defined in (\ref{S epst}) and (\ref{S epss})
are called {\em target} and {\em source counital maps}
and are denoted $\eps_t$ and $\eps_s$ respectively :
\begin{equation}
\eps_t(h) = \epsilon(1\1 h)1\2, \qquad
\eps_s(h) = 1\1 \epsilon(h1\2),
\end{equation}
for all $h\in H$. The images of the counital maps
\begin{equation}
H_t = \eps_t(H),  \qquad H_s = \eps_s(H)
\end{equation}
are separable subalgebras of $H$, called target and source {\em bases} 
or {\em counital subalgebras} of $H$.  
These subalgebras commute with each other; moreover
\begin{eqnarray*}
H_t &=& \{(\phi\otimes \id)\Delta(1) \mid \phi\in H^* \} 
= \{ h\in H \mid \Delta(h) = \Delta(1)(h\otimes 1) \}, \\
H_s &=& \{(\id \otimes \phi)\Delta(1) \mid \phi\in H^* \}
= \{ h\in H \mid \Delta(h) = (1\otimes h)\Delta(1) \},
\end{eqnarray*}
i.e., $H_t$ (respectively, $H_s$) is generated 
by the right (respectively, left) tensor factors of $\Delta(1)$
in the shortest possible representation of $\Delta(1)$ in $H\otimes H$.

\subsection*{Integrals and Frobenius weak Hopf algebras}
The notion of an integral in a weak Hopf algebra is a generalization
of that of an integral in a usual Hopf algebra \cite[2.1.1]{M1}.

\begin{definition}[\cite{BNSz}]
\label{integral}
A left (right) {\em integral} in $H$ is an element
$\ell\in H$ ($r\in H$) such that 
\begin{equation}
h\ell =\eps_t(h)\ell, \qquad (rh = r\eps_s(h)) \qquad \mbox{ for all } h\in H. 
\end{equation}
\end{definition}

We denote by $\int_H^l$ (respectively, by
$\int_H^r$) the space of left (right) integrals in $H$.  
Clearly, $\int_H^l$ is a left ideal of $H$ and $\int_H^r$ is a right
ideal of $H$.

An integral in a finite dimensional weak Hopf algebra
$H$ (left or right) is called {\em non-degenerate} if
it defines a non-degenerate functional on $H^*$. A left integral $\ell$
is called {\em normalized} if $\eps_t(\ell)=1$. Similarly, $r\in \int_H^r$ 
is normalized if $\eps_s(r)=1$. 

Any left integral $\lambda\in \int_{H^*}^l$ satisfies the following
invariance property :
\begin{equation}
\label{left invariance}
g\1 \la \lambda,\, hg\2 \ra = S(h\1) \la \lambda,\, h\2g \ra,
\qquad g,h\in H.
\end{equation}
Similarly, a right integral $\rho\in \int_{H^*}^r$ satisfies
\begin{equation}
\label{right invariance}
\la \rho,\, g\1h \ra  g\2  = \la \rho,\, gh\1 \ra S(h\2),
\qquad g,h\in H.
\end{equation}

The fundamental theorem for weak Hopf modules \cite[3.9]{BNSz}
implies that if $H$
is finite dimensional then $\int_{H}^l$ generates $H$ as a dual left
$H^*$-module. More precisely, there is an isomorphism $H\cong
\int_{H}^l \otimes_{H_t^*} H^*$ of right $H^*$-Hopf modules given by
\begin{equation*}
\ell \otimes \phi \mapsto S(\phi) \actl \ell, 
\end{equation*}
where $\ell\in \int_{H}^l$ and
$\phi\in H^*$. Note that, unlike in the usual Hopf algebra case, this
does not immediately imply existence of non-degenerate integrals in $H$. 
However, it follows from the arguments of \cite[Section 3]{BNSz} that $H$ is
automatically a Frobenius algebra when its bases are commutative.
It was also shown in \cite[3.16]{BNSz} that $H$ is a Frobenius algebra
if and only if there exists a non-degenerate integral in $H$,
if and only if $\dim \int_{H}^l =\dim H_t$, and if and only if
$H^*$ is a Frobenius algebra. 

\begin{note}
In what follows we always work with Frobenius weak Hopf algebras.
\end{note}

Maschke's theorem for weak Hopf algebras, proved in \cite[3.13]{BNSz}
states that a weak Hopf algebra $H$ is semisimple if and only if  $H$ is
separable and if and only if  there exists a normalized left integral in $H$.
In particular, every semisimple weak Hopf algebra is finite dimensional. 
Note that since 
$\eps_t(\int_{H}^l) \subset Z(H) \cap H_t$, it follows that
$H$ is semisimple if and only if there exists $\ell\in \int_{H}^l$
such that $\eps_t(\ell)p \neq 0$ for all primitive idempotents 
$p\in Z(H) \cap H_t$.

For a Frobenius $H$ , there is a useful notion of duality between 
left integrals in $H$ and $H^*$ \cite[3.18]{BNSz}. If $\ell\in \int_{H}^l$
is a left integral and there exists $\lambda\in H^*$ such that
$\lambda\actl \ell = 1$, then such $\lambda$ is unique, it is a left integral
in $H^*$, and  both $\ell$ and $\lambda$ are non-degenerate. Moreover,
$\ell\actl \lambda =\epsilon$. Such a pair of integrals is called a pair
of {\em dual} integrals.
\end{section}

\begin{section}
{Minimal weak Hopf algebras}

Here we classify  all weak Hopf algebras generated by their bases. 
These are probably the most trivial (although in general
non-commutative and non-cocommutative) weak Hopf algebras, as every
weak Hopf algebra contains a weak Hopf subalgebra with this property.

Let $H$ be a weak Hopf algebra.
Observe that a subalgebra $H_{\text{min}} := H_t H_s$ is 
a finite dimensional weak Hopf subalgebra of $H$ (this follows
from the fact that bases are finite dimensional and commute with each other).
Moreover, since $H_t$ and $H_s$ are generated
by the tensor factors  in the shortest possible representation of $\Delta(1)$,
$H_{\text{min}}$ is the minimal weak Hopf subalgebra of $H$ containing $1$.

\begin{definition}
\label{min wha}
A weak Hopf algebra is called {\em minimal} if 
it has no proper weak Hopf subalgebras.
\end{definition}

\begin{remark}
The only minimal usual Hopf algebra is $k1$.
\end{remark}

We have
\begin{equation*}
H_{\text{min}} \cong H_t \otimes_{H_t \cap H_s} H_s \cong
H_t \otimes_{H_t \cap H_s} H_t^{\op},
\end{equation*}
as algebras (with the obvious multiplication in the relative
tensor product), where $H_t \cap H_s \subseteq  Z(H_t)\cap Z(H_s)$.

Recall that for an algebra $A$ an element $\mathcal{E}\in A\otimes A$
is called a {\em separability element} if $(a\otimes 1)\mathcal{E}
= \mathcal{E}(1\otimes a)$ for all $a\in A$ and $m(\mathcal{E})=1$,
where $m$ is the product in $A$. An algebra $A$ for which there exists
such an element is called {\em separable}; it is a standard fact in
the theory of associative algebras that a separable algebra over a
field is finite dimensional and semisimple. A separability element
$\mathcal{E}$ is {\em two-sided} if $\mathcal{E}(a\otimes 1)
= (1\otimes a)\mathcal{E}$; if a  two-sided  separability
element exists, it is unique.

\begin{remark}
Note that $H_{\text{min}}$ is a separable algebra with a separability
element $\mathcal{E} = 1\1 S(1'\1) \otimes S(1\2) 1'\2$.
In particular, any minimal weak Hopf algebra is semisimple. 
\end{remark}

Assume that $k$ is an algebraically closed field of characteristic $0$.
The next Proposition gives an explicit description and complete
classification of minimal weak Hopf algebras over $k$.

\begin{proposition}
\label{claasification of minimals}
Every minimal weak Hopf algebra $H$ is completely determined
by the following data : $(B,\,A,\, g)$, where $B$
is a finite dimensional semisimple algebra, 
$A\subseteq Z(B)$ is a commutative subalgebra, 
and $g\in B$ is an invertible element such that
$\Tr(\pi(g))=\deg\pi$ for any finite dimensional representation 
$\pi$ of $B$, as follows.

As an algebra, $H\cong B\otimes_A B^{\op}$, i.e., $H$  is generated
by the elements $b\in B$, $\ov{c}\in B^{\op}$ (where $c\mapsto \ov{c}$
is the canonical algebra anti-isomorphism from $B$ to $B^{\op}$) and
relations $b \ov{c} = \ov{c} b$, $a =\ov{a},\,a\in A$.
The coalgebra structure is given by
\begin{eqnarray}
\Delta(b \ov{c}) &=& b (\ov{g e^{(1)}}) \otimes e^{(2)} \ov{c},\\
\epsilon(b \ov{c}) &=& \Tr_{\text{reg}}(g^{-1}cb),
\end{eqnarray}
where $e = e^{(1)} \otimes e^{(2)}$  is the unique $2$-sided separability
element of $B$ and $\Tr_{\text{reg}}$ is the trace of the regular 
representation of $B$, and the antipode is given by
\begin{equation}
S(b \ov{c}) = g^{-1}cg \ov{b}.
\end{equation}
Minimal weak Hopf algebras defined by $(B,\,A,\, g)$ and $(B',\,A',\, g')$
are isomorphic if and only if there exists an algebra isomorphism
$\tau :B \to B'$ such that $\tau(A) = A'$ and $\tau(g) =g'$.
\end{proposition}
\begin{proof}
Let $B =H_t$ and $A = H_t \cap H_s \subseteq Z(B)$, then we clearly have
$H \cong  B\otimes_A B^{\op}$ as algebras (where $B^{\op}\cong H_s$).
We identify the map $b\mapsto \ov{b}$ with the antipode map of $H$
restricted to $H_t$.

It can be easily deduced from the axioms of a weak Hopf algebra 
that for all $b,c\in B$ we have $\epsilon(bc) =\epsilon(cS^2(b))$
and that $\epsilon|_B$ is non-degenerate.
Therefore, there exists an invertible element $g\in B$ such that
$\epsilon(b) = \Tr_{\text{reg}}(g^{-1}b)$. Hence, $S^2(b) = g^{-1}bg$
and $S(b \ov{c}) = S^2(c) S(b) = g^{-1}cg \ov{b}$. 
Since $\eps_t|_{H_s} = S|_{H_s}$, we have :
\begin{equation*}
\epsilon(b \ov{c}) 
= \epsilon(bS(c)) = \epsilon(b S^2(c))
= \epsilon(bg^{-1}cg) = \Tr_{\text{reg}}(g^{-1}cb).
\end{equation*}
Since $S(1\1)\otimes 1\2$ is a separability element of $B$, we have
\begin{equation*}
S(1\1)\otimes 1\2 = e^{(1)} \otimes he^{(2)},
\end{equation*}
where $h\in B$ is such that $e^{(1)}he^{(2)} =1$, i.e., 
$\Tr_{\text{reg}}(\pi(h))=\deg\pi$ for any irreducible representation
$\pi$ of $B$. The counit and antipode properties imply that $h=g$, since
\begin{equation*}
1 = 1\1 \epsilon(1\2) = S^{-1}(e^{(1)}) \Tr(g^{-1}he^{(2)}) = S^{-1}(g^{-1}h).
\end{equation*}
We compute :
\begin{eqnarray*}
\Delta(1)
&=& S^{-1}(e^{(1)}) \otimes g e^{(2)} 
= S(S^{-2}(e^{(1)})) \otimes g e^{(2)} \\
&=& S(ge^{(1)} g^{-1}) \otimes g e^{(2)} 
= \ov{ (ge^{(1)} g^{-1})}  \otimes g e^{(2)} = \ov{(ge^{(1)})}\otimes e^{(2)}.
\end{eqnarray*}
The comultiplication on bases is determined by the value of $\Delta(1)$, 
therefore
\begin{equation*}
\Delta(b \ov{c}) = b (\ov{g e^{(1)}}) \otimes e^{(2)} \ov{c}, \quad
b\in B,\, \ov{c}\in B^{\op},
\end{equation*}
which completely defines the structure of $H$. One can check that the
above operations indeed define a weak Hopf algebra.

Let $H'$ be the minimal weak Hopf algebra defined by  $(B',\,A',\, g')$
and $\tau : H \to H'$ be an isomorphism of weak Hopf algebras.
Then 
\begin{equation*}
\tau(B) = B'\qquad \mbox{and} \qquad  
\tau(A) = \tau(H_t \cap H_s) = H'_t \cap H'_s = A'. 
\end{equation*}
Since $\tau$ preserves the regular trace of $H_t$ and $\epsilon$ 
we have $\tau(g)= \tau(g')$.
\end{proof}

\begin{remark}
Examples of minimal weak Hopf algebra structures on the algebra
$B\otimes_k B^{\text{op}}$ appeared in \cite[Appendix]{BNSz}.
\end{remark}

We will denote the minimal weak Hopf algebra defined by the data
$(B,\,A,\, g)$ by $H_{\text{min}}(B,A,g)$.

\begin{remark}
\label{properties of minimals}
\begin{enumerate}
\item[(i)] The square of the antipode of  $H_{\text{min}}(B,A,g)$
is given by the adjoint action of $g^{-1} S(g)$.  In particular, in any 
weak Hopf algebra with commutative bases, $S^2$ is trivial on the
minimal weak Hopf subalgebra.
\item[(ii)]  For any given $A$ and $B$ as above, $H_{\text{min}}(B,A,1)$
is a unique, up to an isomorphism, minimal weak Hopf algebra, for which
$H_t \cong B$, $H_t\cap H_s \cong A$, and $S^2 =\id$.
\item[(iii)]
For any minimal idempotent $p\in A$, the space $pH_{\text{min}}(B,A,g)$
is a simple subcoalgebra of $H_{\text{min}}(B,A,g)$; in particular,
any minimal weak Hopf algebra is cosemisimple. Duals of weak Hopf
algebras $H_{\text{min}}(B,k1,1)$, which are simple as algebras,
were considered in \cite{NV1}.
\end{enumerate}
\end{remark}

\begin{remark}
\label{S2 =id is ok}
Let us explain that
every weak Hopf algebra is a deformation of  a weak Hopf algebra $H$
with the property that 
\begin{equation}
\label{S2 on bases}
S^2 =\mbox{id}\qquad \mbox {on the minimal weak Hopf subalgebra }
H_{\text{min}} \mbox{ of  } H.
\end{equation}
Explicitly, if $H$ is a weak Hopf algebra and $q\in H_t$ is an invertible
element such that $S^2(q)=q$ and $S(1\1)q 1\2 = 1$ then there is a 
deformation weak Hopf algebra $H_q$ with the underlying algebra $H$ 
and the structure operations
\begin{equation*}
\Delta'(h) = \Delta(h)(1\otimes q), \quad \epsilon'(h) = \epsilon(hq^{-1}),
\quad \mbox{and} \quad S'(h) = q^{-1}S(h)q, \quad h\in H.
\end{equation*}
This deformation can be understood in terms of twisting, cf.\ 
\cite[6.1.4]{NV4}.
In particular, if the minimal weak Hopf subalgebra of $H$ is
$H_{\text{min}}(H_t,H_t \cap H_s,g)$, then $H_{g^{-1}}$ has 
property \eqref{S2 on bases}. Thus, problems regarding general
weak Hopf algebras can be translated to problems regarding those
with the regularity property \eqref{S2 on bases}.

Note that if $H$ satisfies \eqref{S2 on bases} then so does $H^*$.
This appears to be a natural property, since it is satisfied
by dynamical quantum groups \cite{EN1} and weak Hopf algebras
arising as symmetries of Jones-von Neumann subfactors \cite{NV2}, \cite{NV3}.
\end{remark}

The next Corollary shows that, in contrast with usual semisimple
Hopf algebras (cf.\ \cite{St}), there can be infinitely (even
uncountably) many non-isomorphic semisimple weak Hopf algebras
with the same algebra structure.

\begin{corollary}
Let $B$ be a non-commutative finite dimensional semisimple algebra,
and $g_1,\, g_2\in B$ be invertible elements with different spectra.
Then  weak Hopf algebras
$H_{\text{min}}(B,k1,g_1)$ and  $H_{\text{min}}(B,k1,g_2)$
are not isomorphic.
\end{corollary}
\begin{proof}
It is clear that no automorphism of $B$ can map $g_1$ to $g_2$.
\end{proof}

\begin{conjecture}
The number of  non-isomorphic weak Hopf algebras $H$ with the 
property~\eqref{S2 on bases} is finite in any given dimension.
\end{conjecture}

\end{section}


\begin{section}
{Group-like elements in a weak Hopf algebra}

\subsection*{Definition and properties}
Let $H$ be a weak Hopf algebra. The notion of a group-like element
in a weak Hopf algebra was introduced in \cite{BNSz}.

\begin{definition}
An element $g\in H$ is said to be {\em group-like} if it is invertible
and satisfies
\begin{equation}
\label{grp-like equation}
\Delta(g) = (g\otimes g)\Delta(1) \qquad \mbox{and} \qquad
\Delta(g) = \Delta(1)(g\otimes g).
\end{equation}
\end{definition}

Group-like elements of $H$ form a group which we will denote by $G(H)$.

\begin{lemma}
\label{eps on grrouplikes}
For any $g\in G(H)$ we have $\eps_t(g)=\eps_s(g)=1$ and the element
$S(g)=g^{-1}$ is group-like.
\end{lemma}
\begin{proof}
The identities for counital maps follow from applying 
$\epsilon$ to \eqref{grp-like equation}.
That $S(g)=g^{-1}$ follows from the uniqueness of the inverse element 
and  antipode.
\end{proof}

\begin{lemma}
\label{dual grp-like element}
If $H$ is finite dimensional, then $\gamma\in G(H^*)$ if and only if
it is invertible and satisfies the following two conditions :
\begin{eqnarray*}
\label{dual grp-like}
\label{dual grplike 1}
\la\gamma,\, hg \ra &=& \la \gamma,\, h 1\1 \ra \la \gamma,\,S(1\2)g \ra, \\
\label{dual grplike 2}
\la\gamma,\, hg \ra &=& \la \gamma,\, h S(1\1)\ra \la  \gamma,\,1\2g \ra.
\end{eqnarray*}
for all $h,g\in H$.
\end{lemma}
\begin{proof}
This is a  straightforward dualization of \eqref{grp-like equation}.
Equation \eqref{dual grplike 1} is equivalent to $\Delta(\gamma)
= (\gamma\otimes\gamma)\Delta(\epsilon)$ and equation \eqref{dual grplike 2}
is equivalent to $\Delta(\gamma) = \Delta(\epsilon) (\gamma\otimes\gamma)$.
\end{proof}

\begin{remark}
If $H$ is an ordinary Hopf algebra, then the notion of a group-like
element coincides with the usual one (cf.\ \cite{M1}, \cite{S}).
\end{remark}

\begin{proposition}
\label{grp-likes in Hmin}
Any group-like element in a minimal weak Hopf algebra $H_{\text{min}}$
has the form $g= S(y)y^{-1}$, where $y\in H_s$ is such that $S^2(y)=y$.
\end{proposition}
\begin{proof}
Without loss of generality we may assume that $H_t\cap H_s =k1$, as
every minimal weak Hopf algebra is a direct sum of  weak Hopf algebras
with this property.
Let $g = \sum_{i=1}^n y_iz_i$ be a group-like element in  $H_{\text{min}}$,
where $\{y_i\}\subset H_s$ and $\{z_i\}\subset H_t$ are linearly
independent sets. We may assume that each $y_i$ is invertible. Then
\begin{eqnarray*}
\sum_{i=1}^n\, (z_i\otimes y_i)\Delta(1) 
&=& \sum_{ij=1}^n\,  (y_iz_i \otimes y_jz_j)\Delta(1) \\
&=& \sum_{ij=1}^n\, (z_i \otimes y_jz_j S^{-1}(y_i))\Delta(1),
\end{eqnarray*}
whence $y_i = \left(\sum_{ij=1}^n\, y_jz_j \right) S^{-1}(y_i)$
by linear independence. Therefore, $n=1$ and $g = y_1S^{-1}(y_1)^{-1}$.
It is straightforward to check that this $g$ satisfies the second equality 
of \eqref{grp-like equation} if and only if $S^2(y_1) = y_1$.
\end{proof}

\begin{definition}
Elements of the form $g=S(y)y^{-1}$, $y\in H_s$ are called
{\em trivial} group-like elements.
\end{definition}

Group-like elements give rise to weak Hopf algebra automorphisms.

\begin{proposition}
\label{automorphisms}
If $g\in H$ is a group-like element, then the map $h\mapsto ghg^{-1}$, where
$h\in H$, is a weak Hopf algebra automorphism. 
If $\gamma\in H^*$ is a group-like element, then the map 
$h\mapsto (\gamma\actl h
\actr \gamma^{-1})$ is a weak Hopf algebra automorphism.
\end{proposition}
\begin{proof}
A direct computation.
\end{proof}

\begin{proposition}
\label{trivial are trivial}
If $\xi\in H_s^*$ is invertible and $y = (\xi^{-1}\actl 1)\in H_s$ then
\begin{equation}
\label{xi and y}
(S(\xi)\xi^{-1})\actl h \actr (S(\xi)^{-1}\xi) =
S(y)y^{-1}\, h\, S(y)^{-1}y.
\end{equation}
\end{proposition}
\begin{proof}
It follows from the properties of bases of a weak Hopf algebra that
\begin{eqnarray*}
\lefteqn{(S(\xi)\xi^{-1})\actl h \actr (S(\xi)^{-1}\xi) =}\\
&=& (S(\xi)\actl 1)(1\actr S(\xi)^{-1})\, h\, (\xi^{-1}\actl 1)(1 \actr \xi) \\
&=& (\xi^{-1}\actl 1)^{-1} S(\xi^{-1}\actl 1)\, h\, S(\xi^{-1}\actl 1)^{-1}
    (\xi^{-1}\actl 1),
\end{eqnarray*}
for all $\xi \in H_s^*$.
\end{proof}

\begin{definition}
\label{trivial automorphisms}
We will call a weak Hopf algebra automorphism of $H$ {\em trivial} if
it is defined by (\ref{xi and y}).
\end{definition}

Trivial weak Hopf algebra automorphisms of $H$ form a normal subgroup
$\text{Aut}_0(H)$ of the group $\text{Aut}(H)$ of all weak Hopf algebra 
automorphisms of $H$. Let $\widetilde{\text{Aut}(H)}$ denote the corresponding
quotient group.

\subsection*{The group of group-like elements}

Let us denote 
\begin{eqnarray}
\label{half grplikes}
G_1(H) &:=& \{g \mid \Delta(g) = (g\otimes g)\Delta(1) \}\backslash \{0\}, \\
G_2(H) &:=& \{g \mid \Delta(g) = \Delta(1) (g\otimes g)\}\backslash \{0\},
\end{eqnarray}
then the group $G(H)$ of group-like elements of $H$
consists of invertible elements of $G_1(H)\cap G_2(H)$.

The  group $G_0(H)$ of all trivial group-like elements is a normal
subgroup in $G(H)$. Define $\G(H)=G(H)/G_0(H)$, the quotient group of 
$G(H)$ by $G_0(H)$. Let $g\mapsto \tilde{g}$ denote the canonical projection
from $G(H)$ to $\G(H)$.

It turns out that $\G(H)$ plays more important role
than $G(H)$, as it possesses properties extending those of the group of 
group-like elements of a usual Hopf algebra.

\begin{remark}
We have $G_0(H) = \{ 1\}$ if and only if the bases of $H$ coincide.
\end{remark}

\begin{remark}
\label{isomorphic coalgebras}
For any $g\in G(H)$ the map $x\mapsto gx$ is an isomorphism between
cosemisimple coalgebras $H_{\text{min}}=H_tH_s$ and $H_g:=gH_{\text{min}}$. 
For $g,h\in G(H)$ we have $H_g = H_h$ if and only if $\tilde{g}=\tilde{h}$
in $\G(H)$.
\end{remark}

\begin{corollary}
\label{finiteness of G}
If $H$ is finite dimensional then $\G(H)$ is finite.
\end{corollary}
\begin{proof}
By the previous remark, to every $g\in G(H)$ there corresponds a
cosemisimple subcoalgebra $H_g$ of $\text{Corad}(H)$, the coradical of $H$;
and $H_g = H_h$ if and only if $g$ and $h$ define the same coset. But there
are only finitely many mutually non-equal cosemisimple subcoalgebras of
$\text{Corad}(H)$, so the quotient group is necessarily finite.
\end{proof}

\subsection*{Modules associated with group-like elements} 

In the case of a usual Hopf algebra $H$,
the group-like elements of $H^*$ are precisely the homomorphisms
from $H$ to the ground field $k$ (i.e., $H$-module structures on $k$). 
In this subsection we present an analogous correspondence 
for weak Hopf algebras.

Let $H$ be a finite dimensional weak Hopf algebra. For every $\gamma\in
G_1(H^*)$ we define a linear map $\eps_s^\gamma: H\to H_s$ by setting
\begin{equation}
\label{twisted counital s}
\eps_s^\gamma(h) := \la \gamma,\,h 1\1\ra S(1\2).
\end{equation}
Similarly, for every $\gamma\in G_2(H^*)$ we define a linear map 
$\eps_t^\gamma: H\to H_t$ by 
\begin{equation}
\label{twisted counital t}
\eps_t^\gamma(h) :=   S(1\1) \la \gamma,\, 1\2 h \ra, \qquad h\in H.
\end{equation}

It follows from Lemma~\ref{eps on grrouplikes} that $\eps_s^\gamma$
and $\eps_t^\gamma$ are projections, i.e., $\eps_s^\gamma \circ \eps_s^\gamma
= \eps_s^\gamma$ and  $\eps_t^\gamma \circ \eps_t^\gamma = \eps_t^\gamma$.
These projections are generalizations of counital maps, since we have
$\eps_t^\epsilon =\eps_t$ and $\eps_s^\epsilon = \eps_s$.

\begin{lemma}
\label{gamma modules}
For every $\gamma\in G_1(H^*)$ the source counital subalgebra $H_s$
becomes a right $H$-module via
\begin{equation}
\label{right H module}
y\cdot h := \eps_s^\gamma(yh),\qquad y\in H_s,\, h\in H.
\end{equation}
Similarly, for every $\gamma\in G_2(H^*)$ the target counital subalgebra $H_s$
becomes a left $H$-module via
\begin{equation}
\label{left H module}
h\cdot z := \eps_t^\gamma(hz),\qquad z\in H_t,\, h\in H.
\end{equation}
\end{lemma}
\begin{proof}
We will only prove the first statement since the proof of the second one
is completely similar. We have
\begin{equation*}
y\cdot 1 = \la \gamma,\, y 1\1 \ra S(1\2) = \epsilon(y 1\1)S(1\2) =y,
\end{equation*}
for all $y\in H_s$, using Lemma~\ref{eps on grrouplikes} and properties of
the counit. Also,  for all $g,h\in H$ we compute, using 
Lemma~\ref{dual grp-like element} :
\begin{eqnarray*}
(y \cdot g)\cdot h
&=& \la \gamma,\, yg1\1 \ra (S(1\2)\cdot h) \\
&=& \la \gamma,\, yg1\1 \ra \la \gamma,\,S(1\2)h 1\1' \ra S(1\2')\\
&=& \la \gamma,\, ygh 1\1' \ra S(1\2') = y\cdot (gh),
\end{eqnarray*}
where $1'$ stands for the second copy of $1$.
\end{proof}

Let us denote $H_s^\gamma$ (respectively, $H_t^\gamma$) the $H$-module
from Lemma~\ref{gamma modules}.

\begin{proposition}
\label{grp-likes vs modules}
The correspondence $\gamma \mapsto H_s^\gamma$ is a bijection 
between $G_1(H^*)$ and the set of all right $H$-module structures
on $H_s$ that restrict to the regular $H_s$-module. Moreover, if 
$\gamma_1,\gamma_2 \in G_1(H^*)$, then $H_s^{\gamma_1} \cong H_s^{\gamma_2}$
if and only if $\gamma_2  = \gamma_1 S(\xi)\xi^{-1}$ for some $\xi\in H_s^*$.

There is a also a similar bijection between $G_2(H^*)$ and the set of all 
right  $H$-module structures on $H_t$ that restrict to the regular 
$H_t$-module.

Analogous statements also hold for left $H$-modules.
\end{proposition}
\begin{proof}
If $H_s$ has a structure of a right $H$-module via $y\otimes h\mapsto
y\cdot h$, such that its restricted $H_s$-module is the regular
$H_s$-module, then the functional 
\begin{equation}
\label{module gives grp-like}
\gamma : h \mapsto \epsilon(1\cdot h), \qquad h\in H
\end{equation}
belongs to $G_1(H^*)$  by Lemma~\ref{dual grp-like element},
since
$$
\la\gamma,\, hg \ra
= \epsilon((1\cdot h)\cdot g) 
= \epsilon((1\cdot h) 1\1) \epsilon(S(1\2) \cdot g)
= \la \gamma,\, h 1\1\ra \la  \gamma,\,S(1\2)g \ra,
$$
for all $h,g\in H$.  Clearly, formulas (\ref{right H module}) and
(\ref{module gives grp-like}) establish a bijective correspondence.

Next, $H_s^{\gamma_1}$ and $H_s^{\gamma_2}$ are isomorphic if and 
only if there is an invertible element $v\in H_s$ such that
\begin{equation*}
v S(1\2) \la\gamma_1 ,\, h 1\1 \ra = S(1\2) \la \gamma_2,\, hv1\1 \ra,
\end{equation*}
for all $h\in H$. This is equivalent to
\begin{equation*}
\gamma_2 = S^{-2}(v) \actl \gamma_1 \actr  v^{-1} =
\gamma_1 ( S^{-1}(v)\actl \epsilon)(\epsilon \actr v^{-1}),
\end{equation*}
i.e., $\gamma_2 = \gamma_1 S(\xi)\xi^{-1} $, where $ \xi =  S^{-1}(v^{-1})
\actl \epsilon \in H^*_s$. Hence, $\tilde\gamma_1 
= \tilde\gamma_2$ in $\G(H^*)$.
Conversely, if $\gamma_2 = \gamma_1 S(\xi)\xi^{-1}$ for some
$\xi \in  H^*_s$, then the map $y\mapsto vy$, where 
$v= S(\xi^{-1})\actl 1 \in H_s$, is an isomorphism between $H_s^{\gamma_1}$ 
and $H_s^{\gamma_2}$.

In a similar way one can show that if $z \otimes h \mapsto z \circ h$
is a right $H$-module structure on $H_t$ that restricts to the regular 
$H_t$-module, then
\begin{equation}
\gamma : h \mapsto \epsilon(1\circ h), \qquad h\in H
\end{equation}
defines an element of $G_2(H)$. This correspondence is also 
bijective and $H$-modules corresponding to $\gamma_1, \gamma_2$
are isomorphic if and only if $\gamma_2 = \gamma_1 S(\zeta)\zeta^{-1}$
for some $\zeta\in H_t^*$.
\end{proof}

\begin{remark}
\label{two halfs make a grp-like}
Suppose that we have both $H$-module structures on $H_s$ and $H_t$
described in Proposition~\ref{grp-likes vs modules} with actions of $H$
denoted by $\cdot$ and $\circ$ respectively. Let $\gamma$ and $\gamma'$
be the corresponding group-like elements :
\begin{equation*}
\la \gamma,\,h\ra = \epsilon(1 \cdot h), \qquad
\la \gamma',\,h\ra = \epsilon(1 \circ h), \quad h\in H. 
\end{equation*}
If $ \epsilon(1 \cdot h) =  \epsilon(1 \circ h)$ for all $h\in H$, then
$\gamma =\gamma'$ is a group-like element in $H^*$. Indeed, 
by Proposition~\ref{grp-likes vs modules} it satisfies
both conditions of \eqref{grp-like equation}. Note that
$\eps_s(\gamma) =\epsilon$ since
\begin{equation*}
\la \eps_s(\gamma),\, h\ra = \epsilon(1\cdot \eps_s(h)) =\epsilon(h), 
\qquad h\in H,
\end{equation*}
and, therefore, $S(\gamma)\gamma =  \eps_s(\gamma) =\epsilon$, i.e.,
$\gamma$ is invertible.
\end{remark}

The next Proposition describes the space of self-intertwiners of $H_s^\gamma$,
$\gamma \in G(H^*)$.

\begin{proposition}
\label{self-intertwiners}
The map $T \mapsto T(1)$ is an isomorphism between the algebras
$\End (H_s^\gamma)_H$ and $Z(H)\cap H_s$.
These algebras are also isomorphic to $H_t^* \cap H_s^*$.
\end{proposition}
\begin{proof}
Since $H_s^\gamma$ restricts to the right regular $H_s$-module, every 
self-intertwiner of $H_s^\gamma$ is of the form $y\mapsto T(1)y$. The
condition that $T\circ h = h\circ T$ in $\End(H_s^\gamma)_H$ for all $h$
is equivalent to
\begin{equation*}
\la \gamma,\, hS^{-2}(T(1))\ra = \la \gamma,\, T(1)h \ra \qquad
\text{for all } h \in H.
\end{equation*}
This, in turn, is equivalent to the identity
\begin{equation*}
T(1)\actl \epsilon = \epsilon \actr T(1).
\end{equation*}
This means that $\epsilon \actr T(1) \in H_t^* \cap H_s^*$. Observe that
$\zeta \mapsto (\zeta \actl 1)$ is an algebra isomorphism between
$H_t^* \cap H_s^*$ and $Z(H)\cap H_s$ that maps $\epsilon \actr T(1)$ to
$T(1)$. Clearly, $T\mapsto T(1)$ is an algebra isomorphism.
\end{proof}

\end{section}

\begin{section}
{The Radford formula for $S^4$}

\subsection*{Distinguished group-like elements}

We use the idea of Radford \cite{R1} to define a canonical coset 
of group-like elements in the quotient group $\G(H) = G(H)/G_0(H)$.
Let $H$ be a finite dimensional Frobenius weak Hopf algebra
such that $S^2|_{H_{\text{min}}} =\id$,
and let $\ell$ be a non-degenerate left integral in $H$.

\begin{lemma}[(\cite{BNSz}, 3.17)]
\label{l is cyclic}
The element $\ell$ is cyclic and separating
for the right $H_s$-module (respectively,  $H_t$-module) $\int_H^l$, where
$H_s$ (respectively, $H_t$) acts by the right multiplication.
\end{lemma}
\begin{proof}
The map $\xi \mapsto (\xi\actl 1)$ is a linear isomorphism between $H_s^*$
and $H_s$. For all non-zero $\xi\in H_s^*$ we have $0\neq \xi\actl \ell =
\ell(\xi\actl 1)$, whence for all $y\in H_s$ we have $\ell y =0$ 
if and only if $y=0$, i.e., $\ell$ is separating. It is cyclic since
$\dim \int_H^l = \dim H_s$. The proof of the statement regarding $H_t$
is similar and uses that $\ell \actr \xi = \ell (1\actr \xi)$ for all
$\xi \in H_s^*$.
\end{proof}

By Lemma~\ref{l is cyclic}, any non-degenerate $\ell \in \int_H^l$
gives rise to a right $H$-module structures on $H_s$ and $H_t$ given by
\begin{equation}
\label{modules via l}
\ell y h = \ell (y \cdot h) \qquad \mbox{and}  \qquad 
\ell z h = \ell (z \circ h)
\end{equation}
for all $y\in H_s$, $z\in H_t$, and $h\in H$. These modules are well-defined
and restrict to the regular $H_s$- and $H_t$-modules.

\begin{lemma}
\label{coherence of modules}
We have $\epsilon(1\cdot h) = \epsilon(1\circ h)$ for all $h\in H$.
\end{lemma}
\begin{proof}
What we need to show is that the equality $\ell y =\ell z$
for $y\in H_s$ and $z\in H_t$ implies that $\epsilon(y) = \epsilon(z)$.

By Lemma~\ref{l is cyclic} there exists a unique well-defined linear map
$T: H_t \to H_s$ such that $\ell T(z) = \ell z$ for all $z\in H_t$.
Clearly, $T$ is an algebra anti-homomorphism and the composition $ST$ 
is an automorphism of $H_t$.
Since $\epsilon|_{H_t} = \Tr_{\text{reg}}$ when $S^2|_{H_{\text{min}}} =\id$
by Proposition~\ref{claasification of minimals}, we have
$ \epsilon(T(z)) = \epsilon(ST(z)) = \epsilon(z)$, because any automorphism
of a semisimple algebra preserves the regular trace.
\end{proof}

\begin{corollary}
The above modules associated to
any non-degenerate $\ell \in \int_H^l$ canonically define a group-like
element $\gamma_\ell \in G(H^*)$.
\end{corollary}
\begin{proof}
Follows from Remark~\ref{two halfs make a grp-like} and 
Lemma~\ref{coherence of modules}.
\end{proof}

Let us denote by $W_\ell$ the right $H$-module structure on $H_s$
corresponding to $\ell$.

\begin{lemma}
\label{class is well defined}
Let $\ell'$ be another non-degenerate left integral in $H$ and 
$\gamma_{l'}\in G(H^*)$ be the corresponding  group-like element
as above. Then $\tilde\gamma_\ell=\tilde\gamma_{\ell'}$ in $\G(H)$, i.e.,
$\gamma_\ell$ and $\gamma_{\ell'}$ define the same coset in $G(H)$.
\end{lemma}
\begin{proof}
Since both $\ell$ and $\ell'$ are non-degenerate, we have $\ell'=\ell v$ 
for some 
invertible $v\in H_s$, and the action $\bullet$ of $H$ on $W_{\ell'}$ 
is related to the  action $\cdot$ of $H$ on $W_\ell$ by 
$y \bullet h = v^{-1}(vy\cdot h)$ for all $y\in H_s$. 
This means that right $H$-modules $W_{\ell'}$ and $W_\ell$ are
isomorphic via $y\mapsto vy$. Therefore, 
by Proposition~\ref{grp-likes vs modules}, the corresponding group-like
elements $\gamma$ and $\gamma'$ define the same coset in $\G(H^*)$.
\end{proof}

\begin{corollary}
\label{distingusished coset}
Let $H$ be a finite dimensional Frobenius weak Hopf algebra.
There is a canonically defined coset  $\tilde\alpha$ in $\tilde{G}(H)$
such that $W_\ell \cong H_s^\alpha$ for any non-degenerate integral
$\ell\in \int_H^l$ and any representative $\alpha$ of  $\tilde\alpha$.
\end{corollary}

\begin{definition}
We will call $\tilde\alpha$ the {\em canonical coset} of group-like elements
and any representative $\alpha$ of this coset a 
{\em canonical group-like element}.
\end{definition}

\subsection*{Dual pairs of non-degenerate integrals}

We continue to assume that $S^2|_{H_{\text{min}}} =\id$.
Recall from Section~\ref{Prelim} that left integrals $\ell \in \int_H^l$
and $\lambda\in \int_{H^*}^l$ are dual to each other if $\lambda\actl \ell =1$,
in which case they are both non-degenerate and also satisfy $\ell \actl
\lambda =\epsilon$ (this was shown in \cite[3.18]{BNSz}). 

It turns out that such pairs are closely related to distinguished group-like
elements in $H$ and $H^*$. In this subsection we study this connection,
following \cite{R1}.

Recall the projections $\eps_s^\gamma$ and $\eps_t^\gamma$ from 
\eqref{twisted counital s} and \eqref{twisted counital t}.  
For all $\gamma\in G(H^*)$ define
\begin{eqnarray*}
L_\gamma &:=& \{ h\in H \mid gh =\eps_t^\gamma(g)h
          \quad \text{ for all } g\in H \}, \\
R_\gamma &:=& \{ h\in H \mid hg =h \eps_s^\gamma(g)
          \quad \text{ for all } g\in H \}. 
\end{eqnarray*}
Also,  for all $g\in G(H)$ define
\begin{eqnarray*}
L_g &:=& \{ \phi\in H^*\mid \psi\phi =\eps_t^g(\psi)\phi
          \quad \text{ for all } \psi\in H^* \}, \\
R_g &:=& \{ \phi\in H^*\mid \phi\psi =\phi \eps_s^g(\psi)
          \quad \text{ for all } \psi\in H^* \}, 
\end{eqnarray*}

Clearly, $L_g, R_g$ are subspaces of $H^*$ and $L_\gamma, R_\gamma$ are
subspaces of $H$. Note that $L_1 = \int_{H^*}^l,\, R_1 = \int_{H^*}^r$
and $L_\epsilon = \int_{H}^l,\, R_\epsilon = \int_{H}^r$.

\begin{remark}
\label{RL meaning of a}
By Corollary~\ref{distingusished coset}, every non-degenerate $\ell \in
\int_{H}^l$ belongs to some $R_\alpha$ and also every non-degenerate 
$r\in  \int_{H}^r$ belongs to $L_\alpha$, 
where $\alpha$ is a distinguished  group-like element in $G(H^*)$.
Similarly, every non-degenerate $\lambda \in \int_{H^*}^l$ belongs
to some $R_a$ and every non-degenerate $\rho \in \int_{H^*}^r$
belongs to $L_a$, where $a$ is a distinguished  
group-like element in $G(H)$.
\end{remark}

\begin{proposition}
\label{shift}
Let $g,h\in G(H)$. The map $\phi \mapsto (h\actl\phi)$ is a linear
isomorphism between $L_g$ and $L_{gh^{-1}}$.
\end{proposition}
\begin{proof}
Take $\phi\in L_g, \psi\in H^*$ and compute
\begin{eqnarray*}
\psi(h\actl\phi)
&=& h \actl ((h^{-1}\actl\psi)\phi) \\
&=& h \actl (\eps_t^g(h^{-1}\actl\psi)\phi)\\
&=& \eps_t^g(h^{-1}\actl\psi)(h \actl \phi) \\
&=& S(\epsilon\1) \la \epsilon\2 \psi\1,\, g\ra \la \psi\2,\,h^{-1}\ra
       (h \actl \phi) \\
&=&  S(\epsilon\1) \la \epsilon\2 \psi,\, g h^{-1}\ra (h \actl \phi)\\
&=& \eps_t^{gh^{-1}}(\psi)(h \actl \phi).
\end{eqnarray*}
Here we used the definition of $\eps_t^g$ and that $\Delta(\zeta)=
\epsilon\1\zeta \otimes \epsilon\2$ for all $\zeta\in H_t^*$. Thus,
$(h\actl \phi)\in L_{gh^{-1}}$, which proves the claim.
\end{proof}

Let $\ell$ be a left non-degenerate integral in $H$ and $\lambda$ be the
dual left integral in $H^*$. Then $S(\ell)$ is a right integral in $H$, i.e.,
$S(\ell)\in R_\epsilon$. By Remark~\ref{RL meaning of a} we also have
$S(\ell)\in L_\alpha$, where $\alpha$ is a distinguished group-like element
in $G(H^*)$. By Proposition~\ref{shift} we have $\alpha^{-1}\actl S(\ell)\in
L_\epsilon$, therefore $S(\ell) = \alpha \actl (\ell y_\ell)$ for some
$y_\ell\in H_s$. Below we show that $y_\ell =1$. 

The following Lemma was proved in \cite[2.2]{R1} for ordinary Hopf algebras.
It is valid for weak Hopf algebras as well, we give a proof for the sake of 
completeness.

\begin{lemma}
\label{S via integrals}
Let $\ell,\,\lambda$ be a dual pair of integrals as above.
Then for all $\phi\in H^*$ we have
\begin{equation}
S(\phi) = (\ell \actr \phi)\actl \lambda.
\end{equation}
\end{lemma}
\begin{proof}
We compute, using the invariance property (\ref{left invariance})
of left integrals :
{\allowdisplaybreaks
\begin{eqnarray*}
(\ell \actr \phi)\actl \lambda 
&=& \la \phi\1,\, \ell\1 \ra \ell\2 \actl \lambda \\
&=& \lambda\1 \la \phi\lambda\2,\,\ell \ra \\
&=& S(\phi\1) \la \phi\2\lambda,\,\ell \ra \\
&=& S(\phi\1) \la \phi\2,\,\lambda\actl \ell \ra\\
&=& S(\phi),
\end{eqnarray*}}
for all $\phi\in H^*$.
\end{proof}

\begin{lemma}
\label{y=1}
We have $y_\ell =1$.
\end{lemma}
\begin{proof}
First, we note that
$S(\ell) =\alpha \actl (ly_\ell) = (y_\ell \actl \alpha)\actl \ell$,
by the definition of $y_\ell$. On the other hand,
\begin{equation*}
S(\ell) = (\lambda \actr \ell) \actl \ell
\end{equation*}
by Lemma~\ref{S via integrals}. Since $\ell$ is non-degenerate, we have
\begin{equation}
\label{alpha lambda l}
y_\ell \actl \alpha = \lambda \actr \ell.
\end{equation}
Applying the target counital map to the both sides of the last identity
we get
\begin{eqnarray*}
\eps_t(\lambda \actr \ell)
&=& \la \lambda\1,\, \ell \ra \eps_t(\lambda\2) 
     = \la \epsilon\1\lambda, \ell \ra \epsilon\2  = \epsilon \\
\eps_t(y_\ell \actl \alpha) 
&=& \eps_t(\alpha(y_\ell \actl \epsilon)) 
     = \alpha S(y_\ell\actl \epsilon)\alpha^{-1},
\end{eqnarray*}
therefore $y_\ell\actl \epsilon =\epsilon$, i.e., $y_\ell =1$. 
\end{proof}

\begin{corollary}
\label{integrals meet group-likes}
For any pair of left integrals $\ell$ and $\lambda$ such that
$\ell \actl\lambda =\epsilon$ and $\lambda \actl \ell =1$ there exist
distinguished group-like elements $\alpha\in G(H^*)$ and $a\in G(H)$
such that
\begin{equation}
\lambda \actr \ell =\alpha \qquad \mbox{and} \qquad
\ell \actr \lambda =a.
\end{equation} 
Furthermore, $S(\ell) = \alpha\actl \ell$ and $S(\lambda)=a \actl \lambda$. 
\end{corollary}
\begin{proof}
The statements concerning $\alpha$ follow from \eqref{alpha lambda l}
in the the proof of Lemma~\ref{y=1}
and those concerning $a$ follow by duality.
\end{proof}

\subsection*{Formula for $\mathbf{S^4}$}
We extend the argument of \cite[Section 3]{R1} to establish
a formula for the fourth power of the antipode in terms of
distinguished group-like elements.
We keep notation of the previous subsection. Let $\alpha$ and $a$
be distinguished group-like elements from
Corollary~\ref{integrals meet group-likes}.

Define linear isomorphisms 
$\ell_L,\ell_R : H^*\to H$ and $\lambda_L, \lambda_R : H\to H^*$ by 
\begin{gather*}
\ell_L(\phi) = \phi \actl \ell, \qquad \lambda_L(h) = h\actl \lambda,\\
\ell_R(\phi) = \ell \actr \phi, \qquad \lambda_R(h) = \lambda\actr h, 
\end{gather*}
for all $\phi\in H^*$ and $h\in H$.

\begin{proposition}
\label{lambda-ell maps}
We have the following relations
{\allowdisplaybreaks
\begin{eqnarray*}
\ell_L \circ \lambda_R(h) &=& S(h), \\
\ell_L \circ \lambda_L(h) &=& S^{-1}(\alpha \actl h),\\
\ell_R \circ \lambda_R(h) &=& S^{-1}(a^{-1}h), \\
\ell_R \circ \lambda_L(h) &=& S((h\actr \alpha)a^{-1}),
\end{eqnarray*}}
for all $h\in H$.
\end{proposition}
\begin{proof}
We establish these relations by a series of direct computations that use
the invariant properties of integrals and 
Corollary~\ref{integrals meet group-likes} :
{\allowdisplaybreaks
\begin{eqnarray*}
\ell_L \circ \lambda_R(h)
&=& \la \lambda\1,\, h\ra \ell_L(\lambda\1) \\
&=& \ell\1 \la \lambda,\,h\ell\2\ra \\
&=& S(h\1) \la \lambda,\, h\2\ell\ra = S(h),\\
\ell_L \circ \lambda_L(h)
&=& \ell\1 \la \lambda,\, \ell\2 h\ra \\
&=& S^{-1}(h\1) \la \lambda,\, \ell h\2\ra \\ 
&=& S^{-1}(h\1)\la \alpha,\,h\2\ra = S^{-1}(\alpha \actl h),\\
\ell_R \circ \lambda_R(h)
&=& \la \lambda, h\ell\1\ra \ell\2 \\
&=& \la a^{-1}\actl S(\lambda),\, h \ell\1\ra \ell\2 \\
&=& \la S(\lambda),\, h\ell\1 a^{-1} \ra (\ell\2 a^{-1})a \\
&=& \la S(\lambda),\, h\1 \ell a^{-1} \ra  S^{-1}(h\2)a \\
&=& \la \lambda,\, h\1\ell \ra S^{-1}(h\2)a \\
&=& S^{-1}(h)a = S^{-1}(a^{-1}h), \\
\ell_R \circ \lambda_L(h)
&=& \la \lambda, \ell\1 h \ra \ell\2 \\
&=& \la a^{-1}\actl  S(\lambda),\, \ell\1 h \ra \ell\2 \\
&=& \la S(\lambda),\, \ell\1 h a^{-1} \ra \ell\2 \\
&=& \la S(\lambda),\,\ell h\1 a^{-1}\ra S(h\2 a^{-1}) \\
&=& \la \lambda,\, \ell h\1\ra S(h\2 a^{-1}) \\
&=& \la \alpha,\, h\1 \ra S(h\2 a^{-1}) = S((h\actr \alpha)a^{-1}),
\end{eqnarray*}}
for all $h\in H$.
\end{proof}

\begin{theorem}
\label{S4 formula}
For all $h\in H$ we have
\begin{equation}
S^4(h) = a^{-1}(\alpha\actl h \actr \alpha^{-1})a.
\end{equation}
\end{theorem}
\begin{proof}
From Proposition~\ref{lambda-ell maps} we deduce
\begin{eqnarray*}
\ell_L \circ \lambda_L(S^2(h))
&=& S^{-1}(\alpha\actl (S^2(h)) \\
&=& S(h\1) \la \alpha,\, S^2(h\2)\ra \\
&=& S(\alpha\actl h)\\
&=& \ell_L \circ \lambda_R(\alpha\actl h), \\
\ell_R \circ \lambda_R(S^2(h))
&=& S^{-1}(a^{-1}S^2(h)) \\
&=& S(a^{-1}h) \\
&=& S((g\actr \alpha)a^{-1}) \\
&=& \ell_R \circ \lambda_L(g),
\end{eqnarray*}
for all $h\in H$, where $g= (a^{-1}ha)\actr \alpha^{-1}$. 
Since the maps $\ell_L$ and $\ell_R$ are bijective, it follows that
\begin{equation}
\lambda_L(S^2(h)) =\lambda_R(\alpha \actl h) 
\qquad \text{ and } \qquad
\lambda_R(S^2(h)) = \lambda_L((a^{-1}ha)\actr \alpha^{-1}).
\end{equation}
Combining these identities we get
\begin{eqnarray*}
\lambda_L(S^4(h)) 
&=& \lambda_R(\alpha \actl S^2(h))\\
&=& \lambda_R(S^2(\alpha \actl h)) \\
&=& \lambda_L((a^{-1}(\alpha \actl h)a)\actr \alpha^{-1}). 
\end{eqnarray*}
Replacing $h$ by $h\actr \alpha$ we get the desired formula.
\end{proof}

\begin{remark}
If $S^2|_{H_{\text{min}}} \neq \id$, then according to 
Remark~\ref{S2 =id is ok} there exists a deformation $H_g$ of $H$
for which the antipode satisfies 
\begin{equation*}
{S'}^2(h) = g^{-1}S(g) S^2(h) g S(g)^{-1}, \qquad h\in H,
\end{equation*}
and ${S'}^2|_{H_{\text{min}}} = \id$, where $g\in H_t$, $S^2(g)=g$. 
\end{remark}

\begin{corollary}
\label{S has finite order}
For any finite dimensional Frobenius weak Hopf algebra
there exists an integer $n\geq 1$ such that $S^{4n}$ is a trivial weak
Hopf algebra automorphism of $H$. In other words, $S^2$ has a finite order
modulo a trivial automorphism.
\end{corollary}
\begin{proof}
By the previous remark we may assume that $S^2|_{H_{\text{min}}} = \id$.
Since any weak Hopf algebra automorphism of $H$ preserves the class
of distinguished group-like elements, one can see that 
weak Hopf algebra automorphisms
\begin{equation}
h \mapsto (a^{-1}(\alpha\actl h \actr \alpha^{-1})a)
\quad \mbox{ and } \quad
h \mapsto (\alpha\actl (a^{-1} h a) \actr \alpha^{-1})
\end{equation}
differ by a trivial automorphism, i.e., $\Ad_{a^{-1}}$ and 
$\Ad_{\alpha}^*$ commute in $\widetilde{\text{Aut}(H)}$.
Since both $a$ and $\alpha$ have finite order modulo trivial group-like
elements,the claim follows.
\end{proof}

\begin{corollary}
\label{S truly finite}
If bases of $H$ coincide, or if bases of $H^*$ coincide
(i.e., $H_t = H_s$, or $H_t^* = H_s^*$), then  $\G(H) = G(H)$ 
and $S$ has finite order in $\text{Aut}(H)$.
\end{corollary}
\begin{proof}
Under the given assumption,
all trivial group-like elements are equal to the identity element $1$. 
\end{proof}

\begin{remark}
\label{biunimodular S4}
In the special case when both $H$ and $H^*$ have non-degenerate $2$-sided
integrals, it was shown in \cite[3.23]{BNSz} that $S^4$ is a trivial
automorphism of $H$. This result also follows from Theorem~\ref{S4 formula}
above since in this case distinguished group-like elements of $H$ and $H^*$
are trivial.
\end{remark}

\end{section}


\begin{section}
{Trace formula and semisimplicity}


Let $k$ be an algebraically closed field of characteristic $0$.

\subsection*{Formula for $\mathbf{\Tr(S^2)}$}

In this Section we derive a weak Hopf algebra
analogue of the Larson-Radford formula
for $\Tr(S^2)$ \cite{LR2}. As in the case of usual Hopf algebras, 
this formula turns  out to be closely related with the semisimplicity of 
the corresponding weak Hopf algebra and its dual.

It follows from \eqref{left invariance} that for a non-degenerate
integral $\ell$ of $H$ with the dual integral $\lambda$, the element
$(\ell\2 \actl \lambda) \otimes S^{-1}(\ell\1) \in H^*\otimes H$ 
is the dual bases tensor. This implies that for every $T\in \End_k H$
we have
\begin{equation}
\label{Tr T}
\Tr(T) = \la \lambda, T(S^{-1}(\ell\1)\ell\2).
\end{equation}
In particular, for $T = S^2$, we get the following 
analogue of \cite[Theorem~2.5(a)]{LR2},
with counits replaced by counital maps.

\begin{proposition}
\label{TR S2}
In any finite dimensional $H$ we have
\begin{equation}
\label{Formula TR S2}
\Tr(S^2) = \la \eps_s(\lambda), \eps_s(\ell)\ra.
\end{equation}
\end{proposition}
\begin{proof}
Follows from \eqref{Tr T} for $T = S^2$.
\end{proof}

\begin{remark}
When $H$ is a finite dimensional Hopf algebra, an immediate consequence
of the above formula is that $H$ and $H^*$ both are semisimple if and only
if $\Tr(S^2) \neq 0$ \cite[Theorem~2.5(b)]{LR2}. Note that the situation
is more complicated for weak Hopf algebras, since one can have $\eps_s(l)=0$
while $\eps_t(\ell) \neq 0$ (e.g., for minimal weak Hopf algebras with
non-commutative bases, see Section $3$).
\end{remark}

Recall from Proposition~\ref{self-intertwiners}
that for any weak Hopf algebra $H$ the algebras $H_s \cap Z(H)$
and $H^*_s \cap H^*_t$ are isomorphic via $y \mapsto (y\actl \epsilon)$.

\begin{definition}
\label{connectedness} 
If the $H_s \cap Z(H)=k1$, or, equivalently, $H^*_s \cap H^*_t =k\epsilon$,
we say that $H$ is {\em connected}.
If both $H$ and $H^*$ are connected, we say that $H$ is {\em biconnected}.
\end{definition}

Below we always assume that $S^2|_{H_{\text{min}}} = \id$.

\begin{proposition}
\label{Tr --> ss}
Let $I$ be the set of primitive idempotents of $H_s \cap Z(H)$. 
If $\Tr(S^2|_{pH})\neq 0$ for all $p\in I$, then $H$ is semisimple.
\end{proposition}
\begin{proof}
Fix a non-degenerate left integral $\ell$ of $H$.
Using the dual bases tensor \eqref{Tr T} we see that
$\Tr(S^2|_{pH})\neq 0$ implies that $y_p = p\eps_s(\ell) \neq 0$.

Let $\alpha$ be a group-like element of $H^*$ such that $S(\ell) =\alpha\actl
\ell$. Such an $\alpha$ exists by Corollary~\ref{integrals meet group-likes}
and depends on $\ell$. Using the definition of integrals we compute :
$$
y_p\ell = S(y_p)\ell = S(\eps_s(\ell))\ell p 
= \eps_t(S(\ell))\ell p= S(\ell) \ell p = S(\ell)y_p,
$$
whence
\begin{equation}
\label{Sy l}
y_p \ell = (\alpha\actl \ell)y_p \qquad \text{for all } p\in I.
\end{equation}
Let $\pi = p\actl 1 \in H_s^* \cap H_t^*$
and let $\xi_p\in \pi H_t^*$ be such that $y_p = \xi_p \actl 1$. For all 
$h\in H$ we have $ yh = \xi_p \actl h$ and $hy = S^{-1}(\xi_p) \actl h$. 
Due to this observation and non-degeneracy of $\ell$, we can rewrite equation 
\eqref{Sy l} as
\begin{equation}
\label{S xi alpha}
\xi_p = S^{-1}(\xi_p)\alpha.
\end{equation} 

By Remark~\ref{properties of minimals}(iii),
$H_{\text{min}}^* =H_t^*H_s^*$ is
a cosemisimple subcoalgebra of $H^*$ and $\pi H_{\text{min}}^*$ is a
simple subcoalgebra of $H_{\text{min}}^*$ (since the element 
$\pi$ is a primitive  idempotent in $H_s^* \cap H_t^*$).

Recall from Remark~\ref{isomorphic coalgebras}
that $i_\alpha : \phi \mapsto \phi\alpha$ is a coalgebra automorphism 
of $H^*$ which preserves $H_{\text{min}}^*$ only if $\alpha$ is trivial. 
This is exactly the case here since the relation \eqref{S xi alpha} 
insures that $i_\alpha$ maps $\pi H_{\text{min}}^*$ to
itself for all primitive idempotents $\pi\in H_s^* \cap H_t^*$.
Thus, one can choose $\ell$ in such a way that $\alpha = \epsilon$,
i.e., one can assume that  $\ell = S(\ell)$ (in which case $\ell$ is a
two-sided integral).   Then \eqref{S xi alpha} implies that $\xi_p \in H_s^* 
\cap H_t^*$, therefore, $y_p \in H_s\cap Z(H)$ is a non-zero multiple of $p$.
We conclude that
$y = \eps_s(\ell) = \sum_p\, y_p$ is an invertible central element
of $H$. Hence, $\ell' = y^{-1} \ell$ is a normalized right integral. 
By Maschke's theorem, $H$ is semisimple.
\end{proof}

\begin{corollary}
\label{dual tr --> ss}
Let $H$ be a finite dimensional weak Hopf algebra.
Suppose that $\Tr(S^2|_{H^*\pi}) \neq 0$ for every primitive idempotent
$\pi \in H_s^* \cap H_t^*$. Then $H$ is semisimple.
\end{corollary}
\begin{proof}
Note that $pH$, where $p\in H_s\cap Z(H)$ is a primitive idempotent,
is naturally identified with the dual space of $H^*\pi$, where 
$\pi = p\actl 1$, so that $\Tr(S^2|_{pH}) = \Tr(S^2|_{H^*\pi})$.
\end{proof}

\begin{corollary}
\label{connected case ss}
Let $H$ be a connected finite dimensional weak Hopf algebra.
If $\Tr(S^2)\neq 0$, then $H$ is semisimple.
\end{corollary}

\begin{corollary}
\label{biconnected case ss}
Let $H$ be a biconnected finite dimensional weak Hopf algebra.
If $\Tr(S^2)\neq 0$, then $H$ and $H^*$ are semisimple.
\end{corollary}

\subsection*{Semisimplicity and cosemisimplicity}

It is a well-known theorem due to Larson and Radford \cite{LR2}
that a semisimple Hopf algebra over $k$ is automatically cosemisimple.
The proof of this theorem uses the finiteness of the order of the antipode
and the formula for $\Tr(S^2)$.

Here we prove a similar result for weak Hopf algebras with coinciding bases.

\begin{lemma}
\label{TR S2 on pHp}
Let $H$ be a semisimple weak Hopf algebra. 
Then $\Tr(S^2|_{pHp}) \neq 0$ for all primitive central idempotents of 
$H_{\text{min}}$.
\end{lemma}
\begin{proof}
Note that since $H$ is semisimple, $S^2$ is an inner algebra automorphism 
of $H$ by \cite[3.22]{BNSz}. For all $p\in H_{\text{min}}$ we have
\begin{equation}
\label{tr S2 as a sum}
\Tr(S^2|_{pHp}) = \Tr(S^2|_{pH_0p}) + \sum_{J}\,\Tr(S^2|_{pJp}),
\end{equation}
where $H_0$ is the minimal weak Hopf quotient of $H$ (viewed as a minimal
ideal of $H$) and the sum is taken over the rest of minimal two-sided ideals 
of $H$ (so that each $pJp$ is a simple algebra). Since  
$p$ is a primitive idempotent
of $H_{\text{min}}$, it follows from Corollary~\ref{S has finite order} that
$S^{n}|_{pHp} = \id$ for some $n>0$, so that $\Tr(S^2|_{pJp}) \geq 0$
by the same reasoning as in \cite[Lemma 3.2]{LR2}.
Since $S^2 = \id$ on $H^*_{\text{min}}$, we have
$S^2|_{H_0} = \id_{H_0}$. Therefore, equation \eqref{tr S2 as a sum}
implies that
\begin{equation}
\Tr(S^2|_{pHp}) > 0
\end{equation}
for all primitive idempotents $p\in H_{\text{min}}$.
\end{proof}

\begin{theorem}
Let $H$ be a semisimple weak Hopf algebra such that $H_t = H_s$.
Then $H^*$ is semisimple.
\end{theorem}
\begin{proof}
By Corollary~\ref{S truly finite}, we have $S^{n} =\id$ for some $n>0$.
Also, the bases are necessarily commutative, so that 
$S^2|_{{H_{\text{min}}}} = \id_{{H_{\text{min}}}}$.

By Lemma~\ref{TR S2 on pHp} we have $\Tr(S^2|_{pHp})\neq 0$ for all
primitive idempotents $p$ of $H_s = H_{\text{min}}$. For any dual pair
$(\ell,\, \lambda)$ of non-degenerate integrals, the element
$$
(\ell\2 \actl \lambda) \otimes p S^{-1}(\ell\2)p
$$ 
is the dual bases tensor for $pHp$, therefore,
\begin{equation*}
0 \neq \Tr(S^2|_{pHp}) = \la \eps_t(\lambda),\, p\eps_s(p\ell) \ra.
\end{equation*}
Let $\rho = (\epsilon \actr p) \in H_t^*$, then $\eps_t(\lambda)\rho \neq 0$.
Since the latter holds for all primitive idempotents $\rho\in H_t^*$, we 
conclude that $\eps_t(\lambda) \in Z(H^*)\cap H_s^*$ is invertible. Therefore,
$\lambda' = \lambda \eps_t(\lambda)^{-1}$ is a normalized left integral and,
hence, $H^*$ is semisimple by Maschke's theorem.
\end{proof}
 

\subsection*{Cosemisimplicity of weak Hopf algebras arising from 
dynamical twisting}

As an application of the above results we prove below that weak
Hopf algebras obtained from dynamical twisting of semisimple Hopf
algebras (such as, e.g., finite group algebras) are semisimple.

A {\em twist} for a weak Hopf algebra $H$ \cite[3.1.1]{EN1}
is a pair $(\Theta,\, \bTheta)$ with
\begin{equation}
\label{Theta and bTheta} 
\Theta\in \Delta(1)(H\otimes H), \quad
\bTheta\in (H\otimes H)\Delta(1), \quad
\mbox { and } \quad
\Theta\bTheta= \Delta(1)
\end{equation}
such that $\Delta_\Theta : H \to H\otimes H$ defined by
\begin{equation*}
\Delta_\Theta(h) = \bTheta\Delta(h)\Theta, \qquad h\in H,
\end{equation*}
is a coassociative comultiplication on the algebra $H$ that makes it
a new weak Hopf algebra, which we will denote by $H_\Theta$.
If we write  $\Theta = \Theta\I \otimes \Theta\II$ and
$\bTheta = \bTheta^{(1)} \otimes \bTheta^{(2)}$, where a summation 
is understood, then the antipode of $H_\Theta$ is given by
\begin{equation}
\label{twisted antipode}
S_\Theta(h) = v^{-1}S(h)v, \qquad h\in H,
\end{equation}
where $v = S(\Theta\I)\Theta\II$ and $v^{-1} = \bTheta^{(1)} S(\bTheta^{(2)})$.
The counital maps of $H_\Theta$ were computed in \cite[3.1.2]{EN1} :
\begin{equation}
\label{twisted t}
{\eps_t}_\Theta(h) = \eps(\Theta\I h)\Theta\II, \qquad
{\eps_s}_\Theta(h) = \bTheta^{(1)}  \eps(h \bTheta^{(2)}),
\end{equation}
for all $h \in H_\Theta$; in particular, the twisting may deform the bases
of a weak Hopf algebra.

Recall the definition of a dynamical twist in a Hopf algebra from \cite{EN1}.
Let $U$ be a Hopf algebra and $A$ be a finite Abelian subgroup of the group
of group-like elements of $H$. Let $A^*$ be the group of characters of $A$,
and let $P_\mu,\, \mu\in A^*$  be the minimal idempotents in $kA$.

\begin{definition}
\label{dynamical twist}
A function $J : A^* \to U \otimes U$ with invertible values and
such that $J(\lambda)\Delta(a) = \Delta(a)J(\lambda)$, for all
$\lambda\in A^*$ and $a\in A$, is called a {\em dynamical twist } 
for $U$ if it satisfies the following functional equations :
\begin{gather}
\label{dynamical equation}
(\Delta\otimes \id)J(\lambda) (J(\lambda+h^{(3)}) \otimes 1)
= (\id \otimes\Delta) J(\lambda) (1\otimes J(\lambda) ), \\
(\eps\otimes \id)J(\lambda) = (\id \otimes \eps)J(\lambda)  = 1.
\end{gather}
Here $J(\lambda+h^{(3)}) =\sum_\mu\, J(\lambda+\mu) \otimes P_\mu
\in U\otimes U \otimes A$.
\end{definition}

We refer the reader to \cite{ES} for an introduction to the theory
of dynamical quantum Yang-Baxter equations and dynamical twists.

Let us consider the tensor product weak Hopf algebra
$H = M_{|A|}(k) \otimes U$, where the matrix algebra $M_{|A|}(k)$
is a weak Hopf algebra (in fact, a groupoid algebra) with the basis
consisting of matrix units $\{ E_{\lambda\mu} \}_{\lambda,\mu\in A^*}$
and with the structure operations
\begin{equation}
\Delta(E_{\lambda\mu}) = E_{\lambda\mu} \otimes E_{\lambda\mu},
\quad \eps(E_{\lambda\mu}) =1, \quad 
S(E_{\lambda\mu}) = E_{\mu\lambda}.
\end{equation}
Then for any dynamical twist $J : A^* \to U \otimes U$ the pair
\begin{gather}
\label{thetas}
\Theta = \sum_{\lambda\mu}\, E_{\lambda \lambda+\mu} J\I(\lambda)
           \otimes E_{\lambda\lambda} J\II(\lambda) P_\mu, \\
\label{thetas1}
\bTheta = \sum_{\lambda\mu}\, E_{\lambda+\mu \lambda} J^{-(1)}(\lambda)
           \otimes E_{\lambda\lambda} P_\mu J^{-(2)}(\lambda), 
\end{gather}
is a twist for $H$ \cite[4.2.4]{EN1}, where $J = J\I \otimes J\II$
and $J^{-1} = J^{-(1)}  \otimes J^{-(2)}$.

By \eqref{twisted t}, the counital maps of 
the corresponding twisted weak Hopf algebra $H_\Theta$ are given by
\begin{equation*}
{\eps_t}_\Theta (E_{\alpha\beta}h) = \epsilon(h) \sum_\lambda\, 
E_{\lambda\lambda} P_{\lambda-\alpha}, \qquad
{\eps_s}_\Theta (E_{\alpha\beta}h) = \epsilon(h)  E_{\beta\beta},
\qquad h \in U.
\end{equation*}

Thus, the bases of $H_\Theta$ are 
\begin{eqnarray*}
(H_\Theta)_t &=& \mbox{span} \{ \sum_\lambda\, 
 E_{\lambda\lambda} P_{\lambda-\alpha} \mid \alpha\in A^* \} \\
(H_\Theta)_s &=& \mbox{span} \{ E_{\beta\beta} \mid
\beta\in A^* \}.
\end{eqnarray*}
Note that the bases do not depend on $J$, and that $H_\Theta$ is
biconnected in the sense of Definition~\ref{connectedness}.

\begin{proposition}
\label{cosemisimplicity of dyn twists}
Let $U$ be a semisimple Hopf algebra, and $J: A^* \to U\otimes U$
be a dynamical twist for $U$. The weak Hopf algebra $H_\Theta$ obtained
by the dynamical twisting of $H = M_{|A|}(k) \otimes U$ is cosemisimple,
i.e., its dual is semisimple.
\end{proposition}
\begin{proof}
By Corollary~\ref{biconnected case ss}, it suffices to show that
$\Tr(S_\Theta^2)\neq 0$ in $\End_k(H_\Theta)$. From \eqref{twisted antipode}
we see that
\begin{equation*}
S_\Theta^2(h) = g^{-1} h g, \qquad \mbox{where } g = S(v)^{-1}v, 
\qquad h\in H_\Theta.
\end{equation*}
Therefore, $\Tr(S_\Theta^2) = \sum_\pi\, \Tr(\pi(g))\Tr(\pi(g^{-1}))$,
where the summation is taken over all irreducible representations of
the semisimple algebra $H$. 
From the formulas for $v$ and $v^{-1}$ and explicit formulas
\eqref{thetas} and \eqref{thetas1} we have :
\begin{eqnarray*}
g &=& \sum_{\lambda\mu}\, E_{\lambda\lambda} P_\mu 
  S(J\I(\lambda)) J\II(\lambda) J^{-(2)}(\lambda) S(J^{-(1)}(\lambda)) P_\mu,\\
g^{-1} &=& \sum_{\lambda\mu}\, E_{\lambda+\mu\lambda+\mu} P_{-\mu} 
  S(J\II(\lambda))J\I(\lambda) J^{-(1)}(\lambda) S(J^{-(2)}(\lambda)) P_{-\mu},
\end{eqnarray*}
whence we can compute
\begin{equation*}
\Tr(\pi(g)) = \sum_{\lambda}\, \Tr(\pi( E_{\lambda\lambda}
 S(J\I(\lambda)) J\II(\lambda) J^{-(2)}(\lambda) S(J^{-(1)}(\lambda))))
=\deg\pi,
\end{equation*}
and, similarly, $\Tr(\pi(g^{-1}))= \deg\pi$ so that $\Tr(S_\Theta^2)
= \dim(H_\Theta)\neq 0$.
\end{proof}

\end{section}


\bibliographystyle{ams-alpha}
  
\end{document}